# Formulae for the determination of the elements of the Eötvos matrix of the Earth's normal gravity field and a relation between normal and actual Gaussian curvature.

G. Manoussakis, D. Delikaraoglou

Dionysos Satellite Observatory and Higher Geodesy Laboratory, Department of Surveying, National Technical University of Athens, tel. 0030 – 210 – 7722693, fax 0030 – 210 – 7722670
gmanous@survey.ntua.gr

**Abstract.** In this paper we form relations for the determination of the elements of the Eötvös matrix of the Earth's normal gravity field. In addition a relation between the Gauss curvature of the normal equipotential surface and the Gauss curvature of the actual equipotential surface both passing through the point $P$ is presented. For this purpose we use a global Cartesian system ($X, Y, Z$) and use the variables $X$, and $Y$ to form a local parameterization a normal equipotential surface to describe its fundamental forms and the plumbline curvature. The first and second order partial derivatives of the normal potential can be determined from suitable matrix transformations between the global Cartesian coordinates and the ellipsoidal coordinates. Due to the symmetry of the field the directions of the local system ($x, y, z$) are principal directions hence the first two diagonal elements of the Eötvös matrix with the measure of the normal gravity vector are sufficient to describe the Gauss curvature of the normal equipotential surface and this aspect gives us the opportunity to insert into the elements of the Eötvös matrix the Gauss curvature.

## Introduction

The Eötvös matrix (along with the Bruns and the Burali – Forti matrices) played an essential role for the study of the Earth's gravity field by appropriate measurements of scalar quantities. With the Bruns matrix we can describe the local variation of the gravity vector expressed in a global Cartesian system ($X,Y,Z$), while the Eötvös matrix is a representation of the Bruns matrix in a local Cartesian system ($x, y, z$). The representation of the Eötvös matrix can be made with invariant quantities such as the curvature and torsion of the geodesic lines which are tangent to the parametric lines at a specific point. For this representation we need a special parameterization of a coordinate patch around a point $P$ of an equipotential surface. This is because of the fact that for the specific form of the Eötvös matrix the parametric lines of the equipotential surface must be at least vertical to each other in a small region around the surface point of interest. Here we form relations which determine the second order partial derivative of the normal potential along the East – West direction, the curvature of the plumbline, the measure of the normal gravity vector and the Gauss curvature of the equipotential surface. For this purpose we will give the necessary transformations between a global Cartesian coordinates ($X,Y,Z$) and the ellipsoidal coordinates ($u, β, λ$). Finally it will be clear that the elements of the normal Eötvös matrix can be determined only with the help of the formula of the normal potential.

## Outline of the method

Suppose that we parameterize locally the equipotential surface at $P$ in $X$ and $Y$ coordinates, where ($X, Y, Z$) is a global Cartesian system such that the $Z$-axis being the Earth's mean axis of rotation, the $X$-axis is the intersection of the equator's plane and the meridian plane of Greenwich and the Y-axis makes the system right – handed. We note that the point P is a point of the Earth's surface with coordinates $(X_P, Y_P, Z_P)$ and normal potential $U_P$. Let

$$\bar{s}^{XY} : \Re^2 \supset D^{XY} \to \Re^3 : (X,Y) \to \bar{s}^{XY}(X,Y) \tag{1.1}$$

be a parameterization of the normal equipotential surface around the point P with vector equation

$$\bar{s}^{XY}(X,Y) = (X, Y, Z(X,Y)) \tag{1.2}$$

A coordinate basis at P is

$$\bar{s}_X^{XY}(P) = \left(1, 0, -\frac{U_X}{U_Z}\right)_P \tag{1.3a}$$

$$\bar{s}_Y^{XY}(P) = \left(1, 0, -\frac{U_Y}{U_Z}\right)_P \tag{1.3b}$$

$$\overline{N}(P) = \frac{\bar{s}_X^{XY}(P) \times \bar{s}_Y^{XY}(P)}{\left|\bar{s}_X^{XY}(P) \times \bar{s}_Y^{XY}(P)\right|} \tag{1.4}$$

The quantities of the first and second fundamental form (Manoussakis et al., 2008) can be expressed with the first and second order partial derivatives of the normal potential

$$U(u,\beta) = \frac{GM}{E}\arctan\frac{E}{u} + \frac{1}{2}\omega^2 a^2 \frac{q}{q_0}\left(\sin^2\beta - \frac{1}{3}\right) + \frac{1}{2}\omega^2(u^2 + E^2)\cos^2\beta \tag{1.5}$$

The partial derivatives of first order at P can be found from

$$\begin{bmatrix} U_X \\ U_Y \\ U_Z \end{bmatrix} = \begin{bmatrix} \frac{u\sqrt{u^2+E^2}}{u^2+E^2\sin^2\beta}\cos\beta\cos\lambda & -\frac{\sqrt{u^2+E^2}}{u^2+E^2\sin^2\beta}\sin\beta\cos\lambda & -\frac{\sin\lambda}{\sqrt{u^2+E^2}\cos\beta} \\ \frac{u\sqrt{u^2+E^2}}{u^2+E^2\sin^2\beta}\cos\beta\sin\lambda & -\frac{\sqrt{u^2+E^2}}{u^2+E^2\sin^2\beta}\sin\beta\sin\lambda & \frac{\cos\lambda}{\sqrt{u^2+E^2}\cos\beta} \\ \frac{u^2+E^2}{u^2+E^2\sin^2\beta}\sin\beta & \frac{u}{u^2+E^2\sin^2\beta}\cos\beta & 0 \end{bmatrix} \begin{bmatrix} \frac{\partial U}{\partial u} \\ \frac{\partial U}{\partial \beta} \\ 0 \end{bmatrix} \tag{1.6}$$

And the second order partial derivatives from the transformation

$$\begin{bmatrix} U_{uu} - U_X \dfrac{E^2}{(u^2+E^2)^{3/2}} \cos\beta\cos\lambda - U_Y \dfrac{E^2}{(u^2+E^2)^{3/2}} \cos\beta\sin\lambda \\ U_{\beta\beta} + U_X (u^2+E^2)^{1/2} \cos\beta\cos\lambda + U_Y (u^2+E^2)^{1/2} \cos\beta\sin\lambda + U_Z u\sin\beta \\ 0 \\ U_{u\beta} + U_X \dfrac{u}{(u^2+E^2)^{1/2}} \sin\beta\cos\lambda + U_Y \dfrac{u}{(u^2+E^2)^{1/2}} \sin\beta\sin\lambda - U_Z \cos\beta \\ 0 \\ 0 \end{bmatrix} = \mathbf{M}_1(u,\beta,\lambda) \begin{bmatrix} U_{XX} \\ U_{YY} \\ U_{ZZ} \\ U_{XY} \\ U_{XZ} \\ U_{YZ} \end{bmatrix} \quad (1.6a)$$

where $M_1$ is a 6 x 6 invertible matrix. The Gauss curvature of the normal equipotential surface at the point $P$ is equal to

$$K_{G,n}(P) = \{[-U_{XX}(U_Z)^2 + 2U_{XZ}U_X U_Z - U_{ZZ}(U_X)^2] \cdot \\ \cdot [(-U_{YY})(U_Z)^2 + 2U_{YZ}U_Y U_Z - U_{ZZ}(U_Y)^2] - \\ - [U_Z(U_Y U_{XZ} + U_X U_{YZ} - U_Z U_{XY}) - U_{ZZ}U_X U_Y]\}_P \dfrac{1}{(U_Z(P))^2 |\overline{\gamma}(P)|^4} \quad (1.7)$$

The angle between the normal gravity vector and the equatorial plane is

$$\Phi_N = \arcsin\dfrac{U_Z(P)}{|\overline{\gamma}(P)|} \quad (1.8)$$

Therefore since the radius of curvature of a parallel circle at the point $P$ is equal to

$$R_p = (X_P^2 + Y_P^2)^{1/2} \quad (1.9)$$

from Meusnier's theorem it is possible to determine the principal curvature along the east – west direction and it is equal to

$$k_1(P) = \cos\left(\arcsin\left(\dfrac{U_Z(P)}{|\overline{\gamma}(P)|}\right)\right)(X_P^2 + Y_P^2)^{-1/2} \quad (1.10)$$

### Formulation of the Eötvös matrix

The non – zero elements of the Eötvos matrix of the normal gravity field at the point $P$ are (Marussi 1985) its diagonal elements $U_{xx}$, $U_{yy}$, $U_{zz}$, and $U_{yz}$ where $(x, y, z)$ is a local astronomical system such that its center is at the point $P$ the $z$ – axis is vertical to the tangent plane at $P$ and pointing outwards, the $x$ – axis is pointing east and $y$ – axis is pointing north. Then the values of the mean curvature and Gauss curvature at the point $P$ of the normal equipotential surface are given by

$$J_{normal}(P) = -(U_{xx}(P) + U_{yy}(P))|\overline{\gamma}(P)|^{-1}, \quad K_{G,n}(P) = U_{xx}(P)U_{yy}(P)|\overline{\gamma}(P)|^{-2} \quad (1.11)$$

In addition for the curvature of the normal plumbline it holds that (Manoussakis, 2008)

$$k_{pl}(P) = \{[U_Y(U_{XZ}U_X + U_{YZ}U_Y) + U_Z(U_{ZZ}U_Y - U_{XY}U_X - U_{YY}U_Y - U_{YZ}U_Z)]^2\big|_P +$$
$$+ [U_Z(U_{XX}U_X + U_{XY}U_Y + U_{XZ}U_Z) - U_X(U_{XZ}U_X + U_{YZ}U_Y + U_{ZZ}U_Z)]^2\big|_P +$$
$$+ [U_X(U_{XY}U_X + U_{YY}U_Y + U_{YZ}U_Z + U_{XX}U_X) - U_Y(U_{XY}U_Y - U_{XZ}U_Z)]^2\big|_P\}^{1/2} \cdot$$
$$\cdot \frac{1}{[(U_X)^2 + (U_Y)^2 + (U_Z)^2]_P^{3/2}}$$

(1.12)

and accordingly in the local system (Moritz – Wellenhof, 2006)

$$k_{pl}(P) = -U_{yz}(P)|\overline{\gamma}(P)|^{-1} \tag{1.13}$$

From (1.10), and (1.13) we have (Maroussi, 1985)

$$U_{xx}(P) = |\overline{\gamma}(P)|\cos\left(\arcsin\left(\frac{U_Z(P)}{|\overline{\gamma}(P)|}\right)\right)(X_P^2 + Y_P^2)^{-1/2} \tag{1.14}$$

Hence the Eötvös matrix for the normal gravity field at the point $P$ can be written as

$$Eötvos_N(P) = \begin{bmatrix} U_{xx} & U_{xy} & U_{xz} \\ U_{yx} & U_{yy} & U_{yz} \\ U_{zx} & U_{zy} & U_{zz} \end{bmatrix} =$$
$$= \begin{bmatrix} U_{xx} & 0 & 0 \\ 0 & (U_{xx})^{-1}|\overline{\gamma}|^2 K_{G,n} & -|\overline{\gamma}|k_{pl} \\ 0 & -|\overline{\gamma}|k_{pl} & 2\omega^2 - U_{xx} - (U_{xx})^{-1}|\overline{\gamma}|^2 K_{G,n} \end{bmatrix}_P$$

(1.15)

**Relation between normal and actual Gaussian curvature.**

Let $W$ be the Earth's gravity potential, and $T$ be the disturbing potential i.e. $W = U + T$. Let $(x_1, y_1, z_1)$ be a local astronomical frame on the actual equipotential surface at the point $P$. The transformation between $(x_1, y_1, z_1)$ and $(x, y, z)$ is given by a matrix of the form

$$\begin{bmatrix} x_1 \\ y_1 \\ z_1 \end{bmatrix} = \begin{bmatrix} \varepsilon_{11} & \varepsilon_{12} & \varepsilon_{13} \\ \varepsilon_{21} & \varepsilon_{22} & \varepsilon_{23} \\ \varepsilon_{31} & \varepsilon_{32} & \varepsilon_{33} \end{bmatrix}\begin{bmatrix} x \\ y \\ z \end{bmatrix}, \quad \varepsilon_{ij} << 1, i \neq j,\; i,j = 1,2,3, \varepsilon_{ii} \cong 1 \tag{1.16}$$

Since $\varepsilon_{ij}$ (for $i$ not equal to $j$) are very small numbers after two successive differentiations of the normal potential $U$ (as $U = U(x_1, y_1, z_1)$) and the disturbing potential $T$ (as $T = T(x_1, y_1, z_1)$) neglecting the second order partial derivatives multiplied by $\varepsilon_{ij}$ and $\varepsilon_{ij}^2$ ($i$ not equal to $j$) we can set

$$T_{x_1 x_1} = \varepsilon_{11}^2 T_{xx} = T_{xx},\; T_{x_1 y_1} = \varepsilon_{11} \varepsilon_{22} T_{xy} = T_{xy},\; T_{y_1 y_1} = \varepsilon_{22}^2 T_{yy} = T_{yy}$$
$$U_{x_1 x_1} = \varepsilon_{11}^2 U_{xx} = U_{xx},\; U_{x_1 y_1} = \varepsilon_{11} \varepsilon_{22} U_{xy} = U_{xy},\; U_{y_1 y_1} = \varepsilon_{22}^2 U_{yy} = U_{yy} \quad (1.17)$$

Therefore we form the equation

$$(U_{xx}(P))^{-1} |\overline{\gamma}(P)|^2 K_{G,n}(P) + T_{yy}(P) = W_{yy}(P) \quad (1.18)$$

Multiplying by $W_{xx}(P)$, subtracting the term $(T_{xy}(P))^2$ from both sides, then after some manipulations we have

$$|\overline{\gamma}(P)|^2 K_{G,n}(P) + \begin{vmatrix} T_{xx} & T_{xy} \\ T_{xy} & T_{yy} \end{vmatrix}_P + \begin{vmatrix} T_{xx} & -T_{yy} \\ U_{xx} & U_{yy} \end{vmatrix}_P = |\overline{g}(P)|^2 K_G(P) \quad (1.19)$$

The above equation is a relation between the Gaussian curvature of the normal equipotential surface with potential $U_P$ and the Gaussian curvature of the actual equipotential surface with potential $W_P$ at the point $P$.

## Conclusions

We outlined a method of determining the elements of the Eötvös matrix of the normal gravity field in local Cartesian system $(x, y, z)$ at a point $P$ of an equipotential surface. We showed that it is possible to determine all necessary quantities of the normal equipotential surface by using a global Cartesian system $(X, Y, Z)$ and several transformations between the first and second order partial derivatives of the normal potential in $(X, Y, Z)$ and $(u, \beta, \lambda)$ coordinates. Due to the symmetry of the normal gravity field the Eötvös matrix contains, the Gauss curvature at the point $P$. The diagonal elements of the Eötvös matrix contain only the second order partial derivative of the normal potential in West – East direction hence we characterize this partial derivative as "*mother second order partial derivative*" for the Eötvös matrix. The second order partial derivative along the east – west direction can be easily determined by the relation (1.14). This relation contains only the magnitude of the normal gravity vector and first order partial derivative of the normal potential $U$ along the $Z$ direction. We mention that it is not possible to form an analogous simple relation for the second order partial derivative of the normal potential $U$ along the north south direction. This is because the principal curvature along this direction is described with a quite complicated function and it is not possible to be related with the curvature of the parallel circle passing through point $P$. Finally we gave a relation between the Gauss curvature of a normal equipotential surface, and an actual equipotential surface, which pass from the point $P$. This relation holds if we make certain assumptions which are described from (1.17). These assumptions occur from the transformation matrix between the two local systems which is "close" to the identity matrix.